\documentclass[12pt]{article}

\usepackage{theorem,amssymb}

\advance \topmargin by -\headheight
\advance \topmargin by -\headsep
\evensidemargin \oddsidemargin
\author{J.-P. Allouche\thanks{supported by MENESR, 
ACI NIM 154 Num\'eration.} \\
CNRS, LRI \\
B\^atiment 490 \\
Universit\'e Paris-Sud \\
F-91405 Orsay Cedex \\
France \\
{\tt allouche@lri.fr} \\
}

\title{A note on univoque self-Sturmian numbers}

\date{ }

\def \proof{\bigbreak\noindent{\it Proof.\ \ }}

\def \endpf{{\ \ $\Box$ \medbreak}}

\newtheorem{proposition}{Proposition}
\theorembodyfont{\rm}

\newtheorem{remark}{Remark}

\begin{document}

\maketitle

\begin{abstract}
We compare two sets of (infinite) binary sequences whose suffixes satisfy 
extremal conditions: one occurs when studying iterations of a unimodal 
continuous map from the unit interval into itself, but it also characterizes
univoque real numbers; the other is an equivalent definition of characteristic 
Sturmian sequences. As a corollary to our study we obtain that a real number 
$\beta$ in $(1,2)$ is univoque and self-Sturmian if and only if the 
$\beta$-expansion of $1$ is of the form $1v$, where $v$ is a characteristic 
Sturmian sequence beginning itself in $1$.

\medskip


\medskip


\end{abstract}

\section{Introduction}

The kneading sequences of a unimodal continuous map $f$ from $[0,1]$ into 
itself, with $f(1)=0$ and $\sup f = 1$ are classically studied by first
looking at the combinatorial properties of the kneading sequence of $1$.
Cosnard proved that, using a simple bijection on binary sequences (namely 
mapping the sequence $(x_n)_{n \geq 0}$ to $(y_n)_{n \geq 0}$, where 
$y_n := \sum_{0 \leq j \leq n} x_j \bmod 2$), the set of kneading sequences
of $1$ for all maps $f$ as above, maps to the set $\Gamma$ defined by
$$
\Gamma := \{u = (u_n)_{n \geq 0} \in \{0, 1\}^{\mathbb N}, \ 
\forall k \geq 0, \ \overline{u} \leq S^k u \leq u \}
$$
where $\overline{u}=(\overline{u}_n)_{n \geq 0}$ is the sequence defined by 
$\overline{u}_n := 1-u_n$, where $S^k$ is the $k$th iterate of the shift 
(i.e., $S^k((u_n)_{n \geq 0}) := (u_{n+k})_{n \geq 0}$), and where $\leq$ is the 
lexicographical order on sequences induced by $0 < 1$. (See \cite{Cos85, AllCos83}; 
actually the relevant set there is $\Gamma \setminus \{(10)^{\infty}\}$. See also
\cite{All83} for a detailed combinatorial study of the set $\Gamma$.)

\bigskip

A slight modification of the set $\Gamma$ describes the expansions of
$1$ in bases $\beta$, where $\beta \in (1, 2)$ is a {\em univoque\,} number,
i.e., a number such that $1$ admits only one expansion in base $\beta$
(see \cite[Remark 1, page 379]{ErdJooKom}, see also \cite{AllCos01} 
and the bibliography therein): the set of all the expansions of $1$ in
bases $\beta \in (1,2)$, where $\beta$ is univoque, is the set
$$
\Gamma_1 := \{u = (u_n)_{n \geq 0} \in \{0, 1\}^{\mathbb N}, \ 
\forall k \geq 0, \ \overline{u} < S^k u < u \}.
$$

\begin{remark}
Note that a binary sequence belongs to $\Gamma_1$ if and only if it belongs
to $\Gamma$ and is not purely periodic.
\end{remark}

\bigskip

Other sequences can be defined by extremal properties of their suffixes:
characteristic Sturmian sequences and Sturmian sequences. More precisely
the following results can be found in several papers (see in particular
\cite{Vee1, Vee2, BulSen, Pir, BugDub}; see also the survey \cite{AllGle} 
and the discussion therein).

\medskip

{\it A binary sequence $u = (u_n)_{n \geq 0}$ is characteristic
Sturmian if and only if it is not periodic and belongs to the set $\Xi$,
where 
$$
\Xi := \{u = (u_n)_{n \geq 0} \in \{0, 1\}^{\mathbb N}, \
\forall k \geq 0, \ 0u \leq S^k u \leq 1u \}.
$$
}

{\it A binary sequence $u = (u_n)_{n \geq 0}$ is Sturmian if and only if it 
is not periodic and there exists a binary sequence $v = (v_n)_{n \geq 0}$,
such that $u$ belongs to $\Xi_v$, where
$$
\Xi_v := \{u = (u_n)_{n \geq 0} \in \{0, 1\}^{\mathbb N}, \
\forall k \geq 0, \ 0v \leq S^k u \leq 1v \}.
$$
The sequence $v$ has the property that $1v = \sup_k S^k u$ and 
$0v = \inf_k S^k u$. This is the characteristic Sturmian sequence having
the same slope as $u$.}
  
\begin{remark}
The reader can find everything on Sturmian sequences in \cite[Chapter 2]{Lot}.
A hint for the proof of the two assertions above is that a sequence is 
Sturmian if and only if it is not periodic and for any binary (finite) word
$w$, the words $0w0$ and $1w1$ cannot be simultaneously factors of the sequence;
furthermore a sequence $u$ is characteristic Sturmian if and only if $0u$
and $1u$ are both Sturmian.
\end{remark}

\section{Comparing the sets $\Gamma$ and $\Xi$}

The analogy between the definitions of $\Gamma$ and $\Xi$ suggests the
natural question whether any sequence can belong to their intersection.
The disappointing answer is the following proposition.

\begin{proposition}
A sequence $u \in \{0, 1\}^{\mathbb N }$ belongs to $\Gamma \cap \Xi$ if
and only if it is equal to $1^{\infty}$ or there exists $j\geq 1$ such that 
$u = (1^j0)^{\infty}$.
\end{proposition}

\proof
If the sequence $u$ belongs to $\Gamma$, we have in particular 
$u \geq \overline{u}$.
Hence $u = 1w$ for some binary sequence $w$. If $u$ is not equal to
$1^{\infty}$ (which clearly belongs to $\Gamma \cap \Xi$), let us write 
$u = 1^j 0 z$ for some integer $j \geq 1$ and some binary sequence $z$.
Since $u$ belongs to $\Gamma$ we have $S^{j+1} u \leq u$, i.e., $z \leq u$.
Now $u$ belongs to $\Xi$, thus $S^j u \geq 0u$, i.e., $0z \geq 0u$, hence
$z \geq u$. This gives $z = u$. Hence $u = (1^j0)^{\infty}$, which in
turn clearly belongs to $\Gamma \cap \Xi$. \endpf

\bigskip

The next question is whether a Sturmian sequence can belong to $\Gamma$.
The answer is more interesting.

\begin{proposition}\label{main}
A (binary) Sturmian sequence $u$ belongs to $\Gamma$ if and only if
there exists a characteristic Sturmian sequence $v$ such that $v$ begins
in $1$ and $u = 1v$.
\end{proposition}

\proof 

Let us first suppose that the Sturmian sequence $u$ belongs to $\Gamma$.
As above, since $u$ belongs to $\Gamma$, $u$ begins in $1$. Hence $u = 1w$ 
for some binary sequence $w$. The inequalities $S^k u \leq u$ for all 
$k \geq 0$ imply that $\sup_k S^k u = u$ (the inequality $\geq$
is trivial since $S^0 u = u$). This can be written $\sup_k S^k u =
1w$. On the other hand $u$ is Sturmian, hence there exists a binary sequence 
$v$ such that $u$ belongs to $\Xi_v$. We also know that $v$ is such that
$1v = \sup_{k \geq 0} S^k u$. Hence $v=w$. Now $\inf_k S^k u = 0v = 0w$. But 
$S^k u \geq \overline{u}$ for all $k \geq 0$, since $u$ belongs to $\Gamma$.
Hence $0v = 0w \geq \overline{u} = 0 \overline{w}$, thus $w \geq \overline{w}$,
hence $w$ begins in $1$.

\medskip

If, conversely, $u = 1v$ where $v$ is a characteristic Sturmian sequence
(which actually implies that $u$ is Sturmian) beginning in $1$, we first 
note that $0v \leq S^k v \leq 1v$ for all $k \geq 0$. Hence, immediately,
$S^k u \leq 1v=u$ for all $k \geq 0$ (using that $S^{k+1} u = S^k v$ and
that $S^0 u = u = 1v$). And also $S^k u \geq 0v \geq 
0\overline{v} = \overline{u}$ (using furthermore that $v \geq \overline{v}$
since $v$ begins in $1$, and that $S^0 u = u \geq \overline{u}$ since $u$
begins in $1$). \endpf

\begin{remark}
We see in particular that a Sturmian sequence belonging to $\Gamma$ must 
begin in $11$. This is not surprising since the only sequence belonging to 
$\Gamma$ that begins in $10$ is $(10)^{\infty}$. This is a particular case 
of a lemma in \cite{All83}: {\em if a sequence $t$ belonging to $\Gamma$ 
begins with $m \overline{m}$, where $m$ is a (finite) nonempty binary word, 
then $t = (m \overline{m})^\infty$.}
\end{remark}

\section{Univoque self-Sturmian numbers}

Several papers were devoted to univoque numbers having an extra property.
For example: 

\noindent
- the smallest univoque number in $(1, 2)$ is determined in
\cite{KomLor}; it is related to the celebrated Thue-Morse sequence and
was proven transcendental in \cite{AllCos00}; 

\noindent
- a detailed study of univoque Pisot numbers belonging to $(1, 2)$ can 
be found in \cite{AllFroHar}.

\bigskip

Self-Sturmian numbers were introduced in \cite{ChiKwo}: these are the real
numbers $\beta$ such that the greedy $\beta$-expansion of $1$ is a Sturmian
sequence on some two-digit alphabet. It is tempting to ask which univoque 
numbers are self-Sturmian. We restrict the study to the numbers in $(1, 2)$ 
for simplicity.

\begin{proposition}
The real self-Sturmian numbers in $(1, 2)$ that are univoque are exactly
the real numbers $\beta$ such that $1 = \sum_{n \geq 1} \frac{u_n}{\beta^n}$, 
where $u = (u_n)_{n \geq 0}$ is a binary sequence of the form $u = 1v$, with 
$v$ a characteristic Sturmian sequence beginning in $1$.
\end{proposition}

\proof This is a rephrasing of Proposition~\ref{main}. \endpf

\begin{remark} 
- The equality $1 = \sum_{n \geq 1} \frac{u_n}{\beta^n}$, where 
$u = (u_n)_{n \geq 0}$ is a binary sequence, uniquely determines
the real number $\beta$ in $(1, 2)$.

- Self-Sturmian numbers correspond to Sturmian sequences of the
form $u = 1v$, where $v$ is any characteristic Sturmian sequence
(see \cite[Remark p.~399]{ChiKwo}). All self-Sturmian numbers are 
transcendental \cite{ChiKwo}. 
\end{remark}

\section{Acknowledgments} The author wants to thank Amy Glen for 
her comments on a previous version of this note.

\end{document}